\documentclass[11pt]{article}
\usepackage{verbatim}
\usepackage{mathrsfs}
\usepackage{hyperref}
\hypersetup{
	colorlinks=true,
	linkcolor=red,
	filecolor=magenta,      
	urlcolor=cyan,
}
\usepackage[abs]{overpic}
\usepackage{subcaption}
\RequirePackage{amsmath,amsfonts,amssymb}
\RequirePackage{rotating}
\RequirePackage{graphicx,xcolor}

\newtheorem{rem}{Remark}

\newcommand{\bF}{\mathbf F}

\newcommand{\bK}{\mathbf K}

\newcommand{\bS}{\mathbf S}

\newcommand{\ba}{\mathbf a}

\newcommand{\bk}{\mathbf k}

\newcommand{\bn}{\mathbf n}
\newcommand{\be}{\mathbf e}

\newcommand{\br}{\mathbf r}

\newcommand{\bu}{\mathbf u}

\newcommand{\bx}{\mathbf x}

\newcommand{\cL}{\mathcal L}

\newcommand{\anna}[2][cyan]{\textcolor{#1}{#2}}

\def\cl {\nonumber \\}
\def\el {\nonumber }

\title{Computational study of numerical flux schemes for mesoscale atmospheric flows in a Finite Volume framework}

\begin{document}


\author{Nicola Clinco$^\star$, Michele Girfoglio$^\star$, Annalisa Quaini$^{\star \star}$, Gianluigi Rozza$^\star$}

\maketitle

\begin{center}
\footnotesize
$\star$ mathLab, Mathematics Area, SISSA, via Bonomea, 265, Trieste, I-34136, Italy \\ $\star \star$ Department of Mathematics, University of Houston, 3551 Cullen Blvd, Houston TX 77204, USA
\end{center}

\begin{abstract}
We develop, and implement in a Finite Volume environment, a
density-based approach for the Euler equations 
written in conservative form using density, momentum, and total energy as variables. 
Under simplifying assumptions, these equations are used to describe 
non-hydrostatic atmospheric flow. 
The well-balancing of the approach is ensured by a local hydrostatic 
reconstruction updated in runtime during the simulation 
to keep the numerical error under control.
To approximate the solution of the Riemann problem, we consider
four methods: Roe-Pike, HLLC, AUSM$^+$-up and HLLC-AUSM. 
We assess our density-based approach and compare the accuracy of these
four approximated Riemann solvers using two
two classical benchmarks, namely the smooth rising thermal bubble and the density current. 
\end{abstract}

\section{Introduction}


It is common practice to divide numerical methods for the solution of the Euler equations into two main categories: density-based and pressure-based. 
Both techniques compute
the velocity field from the momentum equation, however they differ in the practical use of  the continuity equation. In density-based approaches, the continuity equation is used to obtain the density field and the pressure field is computed from the equation of state. 
On the other hand, in pressure-based approaches, the pressure field is computed from a pressure correction equation (called Poisson pressure equation), which is obtained by manipulating continuity and momentum equations. See, e.g., \cite{Miettinen2015}.

Traditionally, 
pressure-based solvers have been designed and mostly used for incompressible and weakly compressible flows, while density-based methods were originally developed
for high-speed compressible flows. We are interested in the numerical simulation of non-hydrostatic mesoscale atmospheric flows, which are governed by the mildly compressible Euler equations. 
Hence, in our first attempt to design an efficient solver for these equations, we have chosen 
a pressure-based approach \cite{GQR_OF_clima, GEA_GIR_QUA, CGQR,ArashGirQuaRo,girfoglio2024comparative}.
This solver is part of GEA (Geophysical and Environmental Applications) \cite{GEA}, a Finite Volume-based open-source package specifically designed for the quick assessment of new computational approaches for the simulation of atmospheric and ocean flows.
Despite the fact that atmospheric flows are mildly compressible, the vast majority of the software packages for weather prediction adopts density-based approaches. So, 
although our pressure-based solver yields accurate results when compared to data in the literature, we present here a density-based solver, which is also going to be featured in GEA. 

Density-based approaches typically  
employ a Riemann solver for the numerical approximation of the flux function. 
For many practical applications, it is extremely expensive to solve the exact Riemann problem as it is a fully non-linear system of equations \cite{godunov59}. Thus, approximate Riemann solvers were developed to capture the main features
of the Riemann problem solution at a reduced computational cost. 
See, e.g., \cite{ROE1981357,Roe1985EfficientCA,Toro1994,LIOU1996364,LIOU2006137,HLL83} and references therein. 
In this paper, we compare the results given by four approximate Riemann solvers: Roe-Pike \cite{Roe1985EfficientCA}, HLLC \cite{Toro1994}, AUSM$^+$-up \cite{LIOU199323,LIOU2006137} and HLLC-AUSM  \cite{LIOU199323,LIOU2006137,KitKeiiShima}. 
Since we are interested in  nearly hydrostatic flows, these approximated Riemann solvers are employed within a 
well-balanced scheme, i.e., a scheme that preserves discrete equilibria, inspired from \cite{KAPPELI2014199,bottaKlein2004}.




We test our well-balanced density-based solver against numerical data available
in the literature for two classical benchmarks for mesoscale atmospheric flow: the
smooth rising thermal bubble \cite{restelliPHD2007,ResGir2009} and the density current \cite{ahmadLindeman2007, strakaWilhelmson1993}. We show that when one uses a
relatively coarse grid, there are noticeable differences in the solutions given by the different methods to approximate the flux function, however such differences become less evident as the mesh is refined. 
In \cite{giraldo_2008}, the authors state that there are no discernible differences in the results obtained with different Riemann solvers when using a Spectral Element method or a Discontinuous Galerkin method for the same benchmarks. However, they present only results with a given mesh, presumably the finest considered. 
Since the level of mesh refinement is typically a compromise
between desired accuracy and required computational time 
with the available computational resources, it is important to show how the 
solution changes when different 
Riemann solvers are adopted so that one can make an informed
decision on what Riemann solver to use for a given mesh. We found that, unless the mesh is very fine, the Roe-Pike and HLLC methods give over-diffusive solutions. 
Both the AUSM$^+$-up and the HLLC-AUSM methods are less dissipative
and thus allow for the use of coarser meshes. In particular, the 
HLLC-AUSM method is the one that gives the best comparison with the data
available in the literature, even with coarser meshes.

The rest of the paper is organized as follows. Sec.~\ref{sec:pd} 
presents the model,
i.e., the mildly incompressible Euler equations. Sec.~\ref{sec:space_disc} 
and \ref{sec:time_disc} discuss the space and time discretization, respectively. Numerical results are shown in Sec.~\ref{sec:num_res}, while some concluding remarks are reported in Sec.~\ref{sec:conclusion}.

 %

\section{Problem definition}\label{sec:pd}

We consider the dynamics of dry atmosphere in a fixed spatial domain $\Omega$. Let $\rho$ be the air density, $\bu= (u, v, w)$  
the wind velocity, $p$ the fluid pressure and $e$ the total energy density. Moreover, let $c_v$ be
the specific heat capacity at constant volume, $T$ the absolute temperature, $g$ the gravitational
constant, and $z$ the vertical coordinate. We write the total energy density as
the sum of three contributions: 
\begin{equation}\label{eq:e}
    e = U + K + \Phi, \quad U = c_v T, \quad K = \frac{|\bu|^2}{2}, \quad \Phi = gz,
\end{equation}
where $U$ is the internal energy density, $K$ is the kinetic
energy density, and $\Phi$ is the gravitational energy density.  The compressible Euler equations state the conservation of mass, momentum and energy 
in $\Omega$ over a time interval of interest $(0, t_f]$:
 \begin{align}
    &\frac{\partial \rho}{\partial t} + \nabla \cdot (\rho \bu) = 0   &&\text{in}\,\, \Omega\times (0,t_f], \label{eq:mass}  \\
    &\frac{\partial (\rho \bu)}{\partial t} +  \nabla \cdot (\rho \bu \otimes \bu) + \nabla p   + \rho g \widehat{\bk} = 0  &&\text{in}\,\, \Omega\times (0,t_f], \label{eq:mom} \\
    &\frac{\partial (\rho e)}{\partial t} +  \nabla \cdot (\rho e \bu) + \nabla \cdot (p \bu) 
    = 0 &&\text{in}\,\, \Omega\times (0,t_f],
    \label{eq:EnergySet}
    \end{align}
where $\widehat{\bk}$ is the unit vector aligned with the vertical axis $z$.
Let $\gamma = \frac{c_p}{c_v}$, where $c_p$ is the specific heat capacity at constant pressure. 
System \eqref{eq:mass}–\eqref{eq:EnergySet} is closed by  a thermodynamics equation of state for $p$ which, based on the assumption that dry air behaves like
an ideal gas, is given by 
\begin{equation}
    p = \rho U \left( \gamma -1\right) = \rho c_v T  \left( \gamma -1\right). \label{eq:state}
\end{equation}

For numerical stability, we add an artificial diffusion term to the momentum and energy equations:
\begin{align}
&\frac{\partial (\rho \bu)}{\partial t} +  \nabla \cdot (\rho \bu \otimes \bu) + \nabla p + \rho g \widehat{\bk} - \mu_a \Delta \bu = \boldsymbol{0}  &&\text{in } \Omega \times (0,t_f], \label{eq:momentum_stab}\\
&\frac{\partial (\rho e)}{\partial t} +  \nabla \cdot (\rho \bu e) +   \nabla \cdot (p \bu) - c_p \dfrac{\mu_a}{Pr} \Delta T = 0 &&\text{in } \Omega \times (0,t_f], \label{eq:energy_stab}
\end{align}
where $\mu_a$ is a constant (artificial) diffusivity coefficient and $Pr$ is the Prandtl number. We note that the choice to have $\mu_a$ constant in space and time is for convenience. More sophisticated LES models can be used  (see, e.g., \cite{CGQR, marrasNazarovGiraldo2015}) without
affecting the findings in this article.

We are interested in nearly hydrostatic flows, i.e., flows originating by a small perturbation of the hydrostatic balance condition characterized by pressure $p_0$ and density $\rho_0$ such that:
\begin{equation}
    \nabla p_0 + \rho_0 g \widehat{\bk} = \textbf{0}\label{eq:stevino}
\end{equation}
Thus, we split pressure and density into mean hydrostatic value and fluctuation over the mean:
\begin{align}
    p(\mathbf{x},t) &= {{p}_{0}}(z) + {p}'(\mathbf{x},t), \label{eq:p_split}\\
    \rho(\mathbf{x},t) &= {\rho_{0}}(z) + {\rho}'(\mathbf{x},t),  \label{eq:rho_split}
\end{align}
where ${\bx} = (x,y,z)$ is a point in the computational domain.
Details on the importance of this splitting can be found in,
e.g., \cite{kellyGiraldo2012}.

For the purpose of rewriting eq.~\eqref{eq:mass}, \eqref{eq:momentum_stab} and \eqref{eq:energy_stab} in compact vector form, let $\mathbf{q}=(\rho,\rho \mathbf{u},\rho e)$ be the solution vector. We define the numerical flux 
\begin{equation}
 \bF^{}(\mathbf{q}) = \begin{pmatrix} \rho \bu \\ \rho \bu \otimes \bu + p \mathbf{I} \\ \left( \rho e + p \right) \bu \\  \end{pmatrix}, \label{eq:FluxVector_inv}
\end{equation}
where $\mathbf{I}$ is the identity matrix, and source term
\begin{equation}
 \bS(\mathbf{q}) = \begin{pmatrix} 0 \\ -\rho g \widehat{\bk} \\ 0 \\  \end{pmatrix}. \label{eq:FluxVector}
\end{equation}
Finally, we introduce diffusion vector $\mathbf{\Gamma} = \left(0,  -\mu_a,  -c_p \dfrac{\mu_a}{Pr}\right)$ and tensor
\begin{equation}
\mathbf{D} = \begin{pmatrix} 
0 \quad 0 \quad 0 \\
0 \quad \bu \quad 0 \\
0 \quad 0 \quad T
\end{pmatrix}.
\label{eq:FluxVector_visc}
\end{equation}
Then, we rewrite eq.~\eqref{eq:mass}, \eqref{eq:momentum_stab} and \eqref{eq:energy_stab} as:
\begin{equation}
    \frac{\partial \mathbf{q}}{\partial t} + \nabla \cdot \mathbf{F(\mathbf{q})} - \mathbf{\Gamma} \Delta \mathbf{D}(\mathbf{q}) = \mathbf{S}(\mathbf{q}).  
    \label{eq:discreteSystem_visc}
\end{equation}

A quantity of interest for atmospheric problems is the potential temperature:
\begin{align}
\theta = \frac{T}{\pi}, \quad \pi = \left( \frac{p}{p_g} \right)^{\frac{R}{c_{p}}}, \label{eq:theta}
\end{align}
where $p_g = 10^5$ Pa is the atmospheric pressure at the ground and $R$ is the specific gas constant of dry air. 
Similarly to the pressure and density, we split $\theta$
into mean hydrostatic value $\theta_0$ and fluctuation over the mean $\theta'$:
\begin{align}
\theta(\bx,t) = \theta_0(z) + \theta'(\bx,t). \el 
\end{align}


\section{Space Discretization}\label{sec:space_disc}

For the space discretization, we adopt the finite volume approach. 
We partition the domain $\Omega$ into cells or control volumes $\Omega_i$, with  $i=1,...,N_c$, where $N_c$ is the total number of cells in the mesh. Let $\mathbf{A}_j = A_j \widehat{\mathbf{n}}_j$ be the surface vector of
each face of the control volume, with $\widehat{\mathbf{n}}_j$ the unit normal vector and with $j = 1, \dots, M$. 

Let us start with the hydrostatic balance. Following the framework introduced in \cite{bottaKlein2004, KAPPELI2014199}, we
determine the local hydrostatic reconstructions \emph{within each control volume} $\Omega_i$ in order to obtain a
well-balanced approximation.  
By combining eq.~\eqref{eq:state}, \eqref{eq:stevino} and \eqref{eq:theta}, we obtain the following equation:
\begin{equation}
      \frac{\partial \rho_0}{\partial z} = - \left [ \frac{\rho_0 g}{\left( \frac{R}{p_g}\right)^{\gamma} }+\gamma \rho_0^{\gamma} \left(\frac{\partial \theta_0}{\partial z} \right) \theta^{\gamma-1}  \right] 
      \frac{1}{\gamma \rho_0^{\gamma-1} \theta_0^\gamma}.\label{eq:drhoz}
\end{equation}
To have the local hydrostatic reconstruction within each control volume $\Omega_i$, we integrate eq.~\eqref{eq:drhoz} over the interval $[z_i - \Delta z_i/2, z_i + \Delta z_i/2]$, where $z_i$ is the vertical coordinate of the cell center of $\Omega_i$ and $\Delta z_i$ is the local vertical mesh size, using the condition $\rho_0 (z = z_i) = \rho_0 (z_i)$, which is cell average density. Notice that, alternatively, $\rho_0 (z_i)$ could be computed by 
a \emph{global} hydrostatic reconstruction used for entire vertical columns in case the mesh is structure vertically (as is usually the case). See, e.g.,  \cite{KlempSkama03,SharLeue02}.  
The approach adopted to solve eq.~\eqref{eq:drhoz} could heavely impact the robustness of the numerical method \cite{bottaKlein2004,KAPPELI2014199}. 
Since we focus on nearly hydrostatic flows with uniform background potential temperature, we can analytically integrate eq.  \eqref{eq:drhoz} to obtain: 
\begin{equation}
\rho_{0,i}(z) = \left( \rho_{0}(z_i)^{\gamma-1} - P \frac{\gamma-1}{\gamma} g (z-z_{i})\right)^{\frac{1}{\gamma-1}} \,\,\,\, \text{in} \,\,\,\,[z_i-\Delta z_i/2, z_i + \Delta z_i/2],
\label{eq:rhoRec} 
\end{equation}
where 
\begin{equation}\label{eq:P}
P=\frac{p_{0}(z_i)}{\rho_{0}(z_i)}^{\gamma}.
\end{equation}
 Note that $\rho_{0,i}(z)$ 
denotes the hydrostatic reconstruction of density within the control volume $\Omega_i$. Eq.~\eqref{eq:P} is equivalent to:
\begin{equation}
    p_{0,i}(z) = P \rho_{0,i}(z)^\gamma \,\,\,\, \text{in} \,\,\,\,[z_i-\Delta z_i/2, z_i + \Delta z_i/2]. \label{eq:pGamma}
\end{equation}




The integral form of eq. \eqref{eq:discreteSystem_visc} on each volume $\Omega_i$ is given by:
\begin{equation*}
    \int_{\Omega_i} \frac{\partial \mathbf{q}}{\partial t} d\Omega + \int_{\Omega_i} \nabla \cdot \mathbf{F}^{} d\Omega - \mathbf{\Gamma} \int_{\Omega_i} \Delta \mathbf{D}^{} d\Omega = \int_{\Omega_i} \bS d\Omega.
\end{equation*}
By applying the Gauss-divergence theorem, the equation above becomes:
\begin{equation}
    \int_{\Omega_i} \frac{\partial \mathbf{q}}{\partial t} d\Omega + \int_{\partial \Omega_i} \mathbf{F}^{} \cdot d\mathbf{A} - \mathbf{\Gamma} \int_{\partial \Omega_i} \nabla \mathbf{D}^{} \cdot d\mathbf{A} = \int_{\Omega_i} \bS d\Omega. \label{eq:discreteSystemII}
\end{equation}

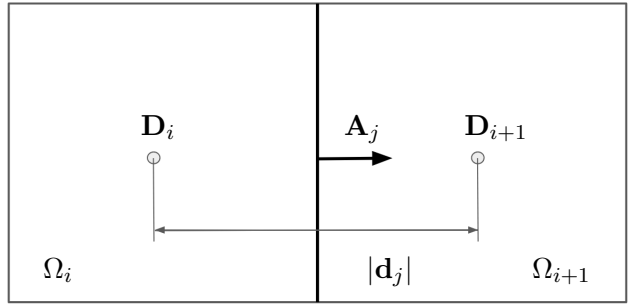
\begin{figure}[htb!]
\centering
\begin{overpic}[abs,unit=1mm,scale=.2, grid=false]{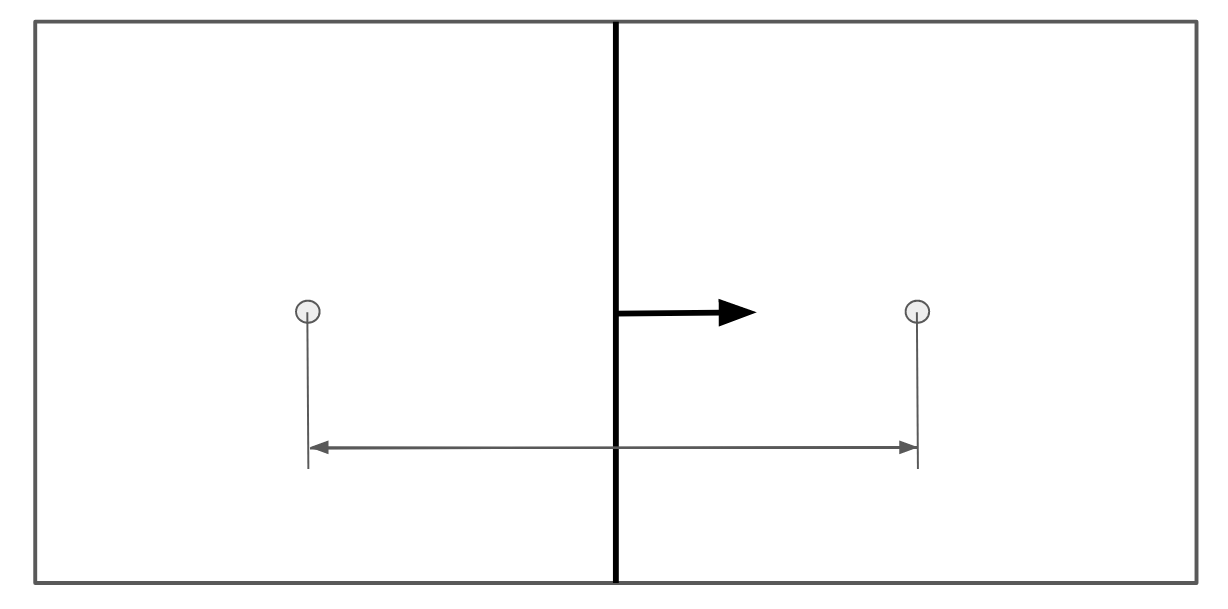}
\put(7,5){\textcolor{black}{{$\Omega_i$}}}
\put(72,5){\textcolor{black}{{$\Omega_{i+1}$}}}
\put(20,24){\textcolor{black}{{$\mathbf{D}_i$}}}
\put(47,24){\textcolor{black}{{$\mathbf{A}_j$}}}
\put(63,24){\textcolor{black}{{$\mathbf{D}_{i+1}$}}}
\put(50,5){\textcolor{black}{{$|\mathbf{d}_{j}|$}}}
\end{overpic}
\caption{Close-up view of two orthogonal control volumes.}
\label{fig:sketch_ortho}
\end{figure}

The discretization of the diffusion term in \eqref{eq:discreteSystemII}, i.e., the
third term at the left hand side, gives:

\begin{equation}
   \int_{\partial \Omega_i} \nabla \mathbf{D}^{} \cdot d\mathbf{A} \approx \sum_j \nabla \mathbf{D}^{}_j \cdot \mathbf{A}_j,
   \label{eq:discreteFlux}
\end{equation}
where $\nabla \mathbf{D}^{}_j$ is the gradient of $\mathbf{D}^{}$ at face $j$. 
On a structured, orthogonal
mesh  (see Fig.~\ref{fig:sketch_ortho}), a second order approximation of $\nabla \mathbf{D}^{}_j$ is given by subtracting the value of $\mathbf{D}$ at the centroid 
of  cell $\Omega_{i}$ from the value of $\mathbf{D}$ at the centroid of
$\Omega_{i+1}$ and dividing by the magnitude of the distance vector $\mathbf{d}_j$ connecting the two
cell centroids:
\begin{equation}\label{eq:diff_term_approx}
\nabla \mathbf{D}^{}_j \cdot \mathbf{A}_j = \dfrac{\mathbf{D}_{i+1} - \mathbf{D}_i}{|\mathbf{d}_j|} |\mathbf{A}_j|.
\end{equation}
For non-structured, non-orthogonal meshes, 
 one has to add an explicit non-orthogonal
correction to the right-hand side of  
\eqref{eq:diff_term_approx}
in order to preserve second order accuracy. See \cite{GQR201927} for details.

For the discretization of the source term in eq.~\eqref{eq:discreteSystemII}, i.e., the term at the right hand side, we follow \cite{bottaKlein2004}. Let 
\begin{equation}
 \bS_0(\mathbf{q}) = \begin{pmatrix} 0 \\ -\rho_0 g \widehat{\bk} \\ 0 \end{pmatrix},
 \quad \widetilde{\bS}_0(\mathbf{q}) = \begin{pmatrix} 0 \\ p_0 \\ 0 \end{pmatrix}.
 \end{equation}
Using \eqref{eq:stevino} and then
applying the Gauss-divergence theorem,
we can write:
\begin{equation}
\int_{\Omega_i} \bS d\Omega \approx \int_{\Omega_i} \bS_0 d\Omega = - \int_{\partial \Omega_i} \widetilde{\bS}_0 \cdot d\textbf{A}.
\approx - \sum_j \widetilde{\bS}_{0,j} \cdot \textbf{A}_j, \label{eq:source}
 \end{equation}
where $\widetilde{\bS}_{0,j}$ is the value of $\widetilde{\bS}_0$ at face $j$ which is computed by using eq. \eqref{eq:pGamma}.

The discretization of the flux term in \eqref{eq:discreteSystemII}, i.e., the
second term at the left hand side, requires more attention. Hence, we treat it in a dedicated subsection.

\subsection{Treatment of the flux term}\label{sec:flux}

\begin{figure}[htb!]
\centering
\begin{overpic}[abs,unit=1mm,scale=.21, grid=false]{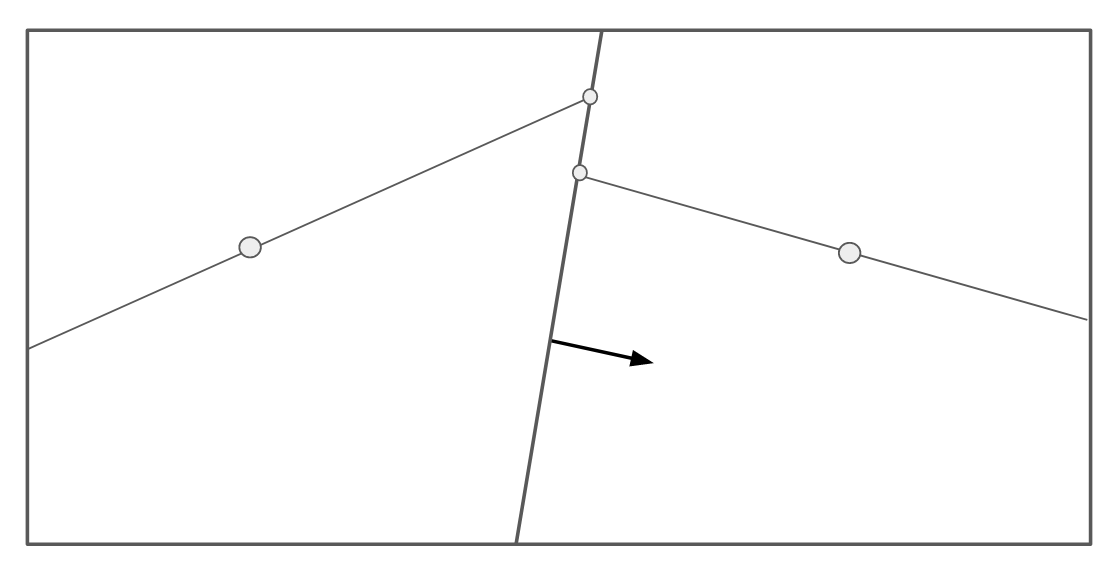}
\put(18,18){{$\mathbf{q}_{i}$}}
\put(45,10){{${\widehat{\bn}}_{j}$}}
\put(35,35){{$\mathbf{q}_{L,j}$}}
\put(44,23){{$\mathbf{q}_{R,j}$}}
\put(61,17){{$\mathbf{q}_{i+1}$}}
\put(8,5){\textcolor{black}{{$\Omega_i$}}}
\put(70,5){\textcolor{black}{{$\Omega_{i+1}$}}}
\end{overpic}
\caption{Close-up view of two non-orthogonal control volumes. 
} 
\label{fig:sketch}
\end{figure}
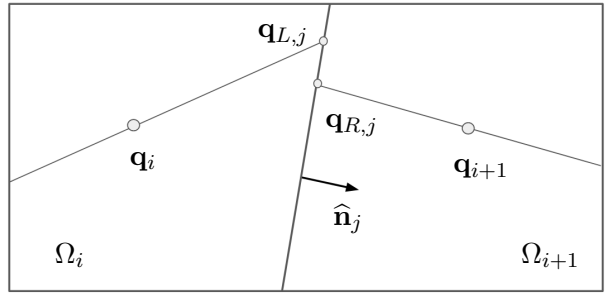

The discretization of the flux term in \eqref{eq:discreteSystemII} gives:

\begin{equation}
   \int_{\partial \Omega_i} \mathbf{F}^{} \cdot d\mathbf{A} \approx \sum_j \mathbf{f}_{j} {A}_j, \label{eq:discreteFlux}
\end{equation}
where $\mathbf{f}^{}_{j} = \mathbf{F}_{j} \cdot \widehat{\mathbf{n}}_{j}$ denotes the numerical flux through face $j$ of $\Omega_i$. We choose to denote $\mathbf{f}^{}_{j}$ as a vector, although for the first component of $\mathbf{F}_{j}$ in \eqref{eq:FluxVector_inv}, i.e., $\rho \bu$, the dot product with $\widehat{\mathbf{n}}_{j}$ gives a scalar, and similarly for the third component of $\mathbf{F}_{j}$. In the case of  incompressible flows, 
the evaluation of $\mathbf{f}^{}_{j}$ requires only interpolation from neighbouring cells. However, 
for compressible flows, fluid properties are not only transported by the flow, but
also by propagation of waves \cite{TadKurg}. 
Thus, the numerical flux is obtained from the solution of a Riemann problem at the cell interfaces
\begin{equation}
\mathbf{f}^{}_{j} = \mathcal{F}(\mathbf{q}_{L,j}, \mathbf{q}_{R,j}),
\label{eq:fluxxx}
\end{equation}
where $\mathcal{F}$ is the adopted Riemann solver and 
$\mathbf{q}_{L,j}$ and $\mathbf{q}_{R,j}$ 
the left and right state at face $j$.
See Fig.~\ref{fig:sketch}, where
$\mathbf{q}_{i}$ and $\mathbf{q}_{i+1}$ are the average solution vector in 
control volumes $\Omega_i$ and $\Omega_{i+1}$. 
Next, we describe how to compute $\mathbf{q}_{L,j}$ and $\mathbf{q}_{R,j}$ to have an equilibrium preserving reconstruction within $\Omega_i$. 


Recall splitting \eqref{eq:p_split} and \eqref{eq:rho_split}
for density and pressure. In each cell $\Omega_i$, these splittings become: 
\begin{align}
    p_i(\mathbf{x},t) &= {{p}_{0,i}}(z) + {p}_{i}'(\mathbf{x},t), \label{eq:pCompI}\\
    \rho_i(\mathbf{x},t) &= {\rho_{0,i}}(z) + {\rho}_{i}'(\mathbf{x},t), \label{eq:rhoCompI}
\end{align}
with ${{p}_{0,i}}(z)$ given in \eqref{eq:pGamma} and 
${\rho_{0,i}}(z)$ in \eqref{eq:rhoRec}. 
The perturbation is reconstructed with
a standard piece-wise linear approximation. 
For the sake of clarity, we present this reconstruction
in 1D. 
Extension to higher dimensions
is not complicated, but notation can become cumbersome. 
In $\Omega_i$, we have:
\begin{align}
    p_i'({x},t) &= p_i'({x}_i,t) + \cL(p'_i) ({x} - {x}_i),
    \label{eq:pCompII}\\
 \rho_i'({x},t) &= \rho_i'({x}_i,t) + \cL(\rho'_i)({x} - {x}_i),
    \label{eq:rhoCompII}
\end{align}
where
\begin{align}
\cL(p'_i) = \mathscr{L}(\mathcal{D}_{i-1}(p'),\mathcal{D}_{i+1}(p')), \label{eq:slope1} \\
\cL(\rho'_i) = \mathscr{L}(\mathcal{D}_{i-1}(\rho'),\mathcal{D}_{i+1}(p')),
\label{eq:slope2}
\end{align}
are the  slopes of the pressure and density perturbation
with the application of a proper limiter $\mathscr{L}$. 
We adopt a monotonized central limiter \cite{VANLEER1979101}. Moreover, following \cite{KAPPELI2014199,bottaKlein2004} we set 
\begin{align}
    \mathcal{D}_{i-1}(p') = \frac{{p}_i'({x}_i) - {p}_{i}'({x}_{i-1})}{{x}_i-{x}_{i-1}}, \quad {p}_{i}'({x}_{i}, t) = p_i({x}_{i}, t) - {{p}_{0,i}}(z_{i}), \label{eq:pert1} \\
    \mathcal{D}_{i-1}(\rho') = \frac{{\rho}_{i}'({x}_{i}) - {\rho}_{i}'({x}_{i-1})}{{x}_i-{x}_{i-1}}, \quad {\rho}_{i}'({x}_{i}, t) = \rho({x}_{i}, t) - {{\rho}_{0,i}}(z_{i}). \label{eq:pert2}
\end{align}

A standard piece-wise linear reconstruction is also applied to the velocity $\mathbf{u}$. 
We omit the formulas for $\mathbf{u}$, since from the formulas above for the reconstruction of a scalar fluctuation 
one can easily
write down the formula for each velocity component. 
Once $\mathbf{u}_i$ is computed in every $\Omega_i$, one uses eq.~\eqref{eq:state} to compute temperature $T_i$
and eq.~\eqref{eq:e} to compute the total energy density $e_i$. 


Finally, we can write the interface solution vectors  $\mathbf{q}_{L,j}$ and $\mathbf{q}_{R,j}$ in 1D:
\begin{align}
    &\mathbf{q}_{L,j}= \mathbf{q}(x_i + \Delta x_i/2), \quad  
    \mathbf{q}_{R,j} = \mathbf{q}(x_{i+1} - \Delta x_{i+1}/2). \el
\end{align}

Now that we have $\mathbf{q}_{L,j}$ and $\mathbf{q}_{R,j}$, 
eq.~\eqref{eq:fluxxx} is completed with the choice of 
 Riemann solver $\mathcal{F}$.
Several options are available in the literature for the computation of $\mathbf{f}^{}_{j}$ \eqref{eq:fluxxx}. 
See, e.g., \cite{ROE1981357,Roe1985EfficientCA,Toro1994,LIOU1996364,LIOU2006137,HLL83}. 
In this paper, we will compare the results given by four approaches: Roe-Pike \cite{Roe1985EfficientCA}, HLLC \cite{Toro1994}, AUSM$^+$-up \cite{LIOU199323,LIOU2006137} and HLLC-AUSM  \cite{LIOU199323,LIOU2006137,KitKeiiShima}.  
These methods are briefly explained in the subsections below.

For simplicity of notation, from now on we will drop the $j$ index from $\mathbf{q}_{L,j}$ and $\mathbf{q}_{R,j}$, i.e., we will use $\mathbf{q}_{L}$ and $\mathbf{q}_{R}$ with the understanding that we are referring to face $j$ in 
cell $\Omega_i$. Similarly, we will drop the 
$j$ index from $\widehat{\mathbf{n}}_{j}$.


\subsubsection{The Roe-Pike method}

Perhaps, the best known of all approximate Riemann solvers is the one due to Roe \cite{ROE1981357}.
The Roe-Pike method represents an improvement over the classical Roe method. See \cite{Roe1985EfficientCA} for details.




In the Roe-Pike method, one first computes the Roe average between left and the right state values:
\begin{align}
    &\widetilde{\rho} = \sqrt{\rho_L \rho_R}, \quad
    \widetilde{\bu} = (\widetilde{u}, \widetilde{v}, \widetilde{w}) = \frac{\sqrt{\rho_L}\bu_L + \sqrt{\rho_R}\bu_R }{\sqrt{\rho_L} + \sqrt{\rho_R}}, \cl
    &\widetilde{h} = \frac{\sqrt{\rho_L}h_L + \sqrt{\rho_R}h_R}{\sqrt{\rho_L} + \sqrt{\rho_R}}, \quad
    \widetilde{a} = \sqrt{(\gamma-1)(\widetilde{h}-K)}, \el
\end{align}
where $h = U + {p}/{\rho} +K $ is the total enthalpy  and $\tilde{a}$ is a modified speed of sound.

Then, one sets 
\begin{align}
    & \widetilde{\alpha}_1 = \frac{1}{2 \widetilde{a}^2} ( p_L -p_R - \widetilde{\rho}\widetilde{a}(u_L -u_R) ), \quad
    \widetilde{\alpha}_2 = (\rho_L -\rho_R) - \frac{p_L - p_R}{\widetilde{a}^2}, \cl
    & \widetilde{\alpha}_3 = \widetilde{\rho}({v}_L - {v}_R), \quad \widetilde{\alpha}_4 = \widetilde{\rho}({w}_L - {w}_R), \quad
    \widetilde{\alpha}_5 = \frac{1}{2 \widetilde{a}^2} ( p_L -p_R + \widetilde{\rho}\widetilde{a}(u_L -u_R) ), \el
\end{align}
and 
\begin{equation}
    \widetilde{\lambda}_1 = \widetilde{u}-\widetilde{a}, \quad \widetilde{\lambda}_2 = \widetilde{\lambda}_3 =   \widetilde{\lambda}_4 = \widetilde{u}, \quad \widetilde{\lambda}_5 =
    \widetilde{u}+\widetilde{a}.\label{eq:RoeEigenvalues}
\end{equation}
With $\widetilde{\ba} = (\widetilde{a}, 0, 0)$, we define 
\begin{equation}
    \widetilde{\br}_1 = \begin{pmatrix} 1 
    \\ \widetilde{\bu}-\widetilde{\ba} \\
    \widetilde{h} - \widetilde{u}\widetilde{a}
    \end{pmatrix},   ~
    \widetilde{\br}_2 = \begin{pmatrix} 
    1 \\ 
    \widetilde{\bu} \\ 
    K
    \end{pmatrix}, ~
    \widetilde{\br}_3 = 
    \begin{pmatrix} 
    0\\ 
    \be_2 \\ 
    \widetilde{v} 
    \end{pmatrix}, ~ 
    \widetilde{\br}_4 =
    \begin{pmatrix} 0 \\ 
    \be_3 \\ 
    \widetilde{w}  
    \end{pmatrix}, ~ 
    \widetilde{\br}_5 = \begin{pmatrix} 
    1 \\
    \widetilde{\bu}+\widetilde{\ba} \\  
    \widetilde{h} + \widetilde{u}\widetilde{a} 
    \end{pmatrix}, \label{eq:RoeEigenVectors}
\end{equation}
where $\be_2$ and $\be_3$ are the second and third
column of the identity matrix of size 3.

At this point, we have all the ingredients
to compute the numerical flux:
\begin{equation}
    \mathbf{f}^{}_{j} = \frac{1}{2}\left(\mathbf{f}^{}(\mathbf{q}_{L}) + \mathbf{f}^{}(\mathbf{q}_{R})\right) -\frac{1}{2}\sum_{k=1}^{5} \widetilde{\alpha}_k{|\widetilde{\lambda}_k|} \widetilde{\br}_k, \label{eq:RoeFlux}
\end{equation}
where
\begin{equation}
\mathbf{f}^{}(\mathbf{q}_{L}) = \begin{pmatrix} \rho_L \bu_L \cdot \widehat{\mathbf{n}} \\ \rho_L \bu_L \otimes \bu_L \cdot \widehat{\mathbf{n}} + p_L  \widehat{\mathbf{n}}\\ \left( \rho_L e_L + p_L\right) \bu_L \cdot \widehat{\mathbf{n}} \\  \end{pmatrix}, \quad 
 \mathbf{f}^{}(\mathbf{q}_{R}) = \begin{pmatrix} \rho_R \bu_R \cdot \widehat{\mathbf{n}} \\ \rho_R \bu_R \otimes \bu_R \cdot \widehat{\mathbf{n}} + p_R \widehat{\mathbf{n}} \\ \left( \rho_R e_R + p_R\right) \bu_R \cdot \widehat{\mathbf{n}} \\  \end{pmatrix}.
 \label{eq:RoeLeftRightFlux}
\end{equation}



\subsubsection{The Harten–Lax–van Leer contact method}


The Harten–Lax–van Leer contact (HLLC) method
is an improvement to the classical Harten–Lax–van Leer Riemann (HLL) solver to solve systems with three or more characteristic fields. It was introduced to avoid the excessive numerical dissipation of HLL for intermediate characteristic fields. We consider a
version of HLLC for the time-dependent Euler equations presented in \cite{Toro1994}.


The definition of $\mathbf{f}_j$  
in this method
requires the left and right wave speed: 
\begin{equation}
    S_L =  \mathbf{u}_L \cdot \widehat{\mathbf{n}} - a_L, \quad S_R =  \mathbf{u}_R \cdot \widehat{\mathbf{n}}  -a_R , \label{eq:onde}
\end{equation}
and the so-called middle wave speed $S^*$ defined as:
\begin{equation}
    S^* = \frac{(p_R - p_L) + \rho_L \mathbf{u}_L \cdot \widehat{\mathbf{n}}(S_L - \mathbf{u}_L \cdot \widehat{\mathbf{n}}) -\rho_R \mathbf{u}_R \cdot \widehat{\mathbf{n}} (S_R-\mathbf{u}_R \cdot \widehat{\mathbf{n}})}{\rho_L(S_L - \mathbf{u}_L \cdot \widehat{\mathbf{n}})-\rho_R(S_R -\mathbf{u}_R \cdot \widehat{\mathbf{n}})}. \label{eq:SstartHLLC}
\end{equation}
In eq.~\eqref{eq:onde}, $a_L$ and $a_R$ are
the left and right speed of sound:
\begin{equation}
    a_L = \sqrt{\gamma \frac{p_L}{\rho_L}}, \quad a_R=\sqrt{\gamma \frac{p_R}{\rho_R}}. \label{eq:vel_sound}
\end{equation}

The numerical flux is computed as: 
\begin{equation}\label{eq:FluxesHLLC}
    \mathbf{f}_{j} = \left\{ 
  \begin{array}{ c l }
    \mathbf{f}(\mathbf{q}_L) & \quad \text{if}  \quad S_L  \geq 0,  \\
    \mathbf{f}(\mathbf{q}^*_L) & \quad \text{if} \quad  S^* \geq 0 \geq S_L,  \\
    \mathbf{f}(\mathbf{q}^*_R) & \quad \text{if} \quad S_R \geq 0 \geq S^*, \\
    \mathbf{f}(\mathbf{q}_R) & \quad \text{if} \quad S_R \leq 0, \\
  \end{array}
  \right.
\end{equation}
where $\mathbf{f}(\mathbf{q}_L)$ and $\mathbf{f}(\mathbf{q}_R)$ are given by eq. \eqref{eq:RoeLeftRightFlux}  and 
\begin{equation}\label{eq:FluxesStarredHLLC}
   \mathbf{f}(\mathbf{q}_\alpha^*) =  \frac{1}{S_L - S^*}\begin{pmatrix}
    (S_\alpha - \mathbf{u}_\alpha \cdot \widehat{\mathbf{n}}) \rho_\alpha \\
    (S_\alpha - \mathbf{u}_\alpha \cdot \widehat{\mathbf{n}}) \rho_\alpha \bu_\alpha +(p^* -p_\alpha) \\
    (S_\alpha - \mathbf{u}_\alpha \cdot \widehat{\mathbf{n}}) \rho_\alpha e_\alpha +(p^*-p_\alpha)S^*\\
  \end{pmatrix},
 \quad \text{with} \quad \alpha=R,L. 
\end{equation}
In \eqref{eq:FluxesStarredHLLC}, $p^*$ is given by  
\begin{equation}
    p^* = \rho_R(\mathbf{u}_R \cdot \widehat{\mathbf{n}}-S_R)(\mathbf{u}_R \cdot \widehat{\mathbf{n}}-S^*) +p_R.
\end{equation}



\subsubsection{The AUSM$^{+}$-up method}

The AUSM$^{+}$-up is an extension of the original advection upstream splitting method (AUSM) aimed at improving the accuracy in applications involving
for low Mach number flows  \cite{LIOU2006137, LIOU1996364}.

The numerical flux $\mathbf{f}_j$ is computed as follows:
\begin{equation}
    \mathbf{f}_{j} = \frac{\dot{m} + |\dot{m}|}{2} \mathbf{\Phi_L} + \frac{\dot{m} - |\dot{m}|}{2} \mathbf{\Phi_R} + \widetilde{p} \mathbf{N}, \label{eq:fluxAUSM+-up}
\end{equation}
with
\begin{equation*}
 \mathbf{\Phi_L} = \begin{pmatrix} 
 1 \\ \bu_L \\ h_L \\  \end{pmatrix}, \quad
 \mathbf{\Phi_R} = \begin{pmatrix} 1 \\ \bu_R \\ h_R \\  \end{pmatrix}, \quad
 \mathbf{N} = \begin{pmatrix} 0 \\ \widehat{\mathbf{n}} \\ 0 \\  
 \end{pmatrix}.
\end{equation*}
In eq.~\eqref{eq:fluxAUSM+-up} $\dot{m}$ is the mass flow rate and $\widetilde{p}$ is the so-called interface pressure \cite{LIOU2006137, LIOU1996364}. 
Below, we explain how to compute these two quantities. 

Let us start with $\widetilde{p}$.
Let $M_L$ and $M_R$ be the normal left and right Mach number
\begin{equation}
 M_L = \dfrac{\mathbf{u}_L \cdot \mathbf{\widehat{n}}}{a_{1/2}}, \quad M_R = \dfrac{\mathbf{u}_R \cdot \mathbf{\widehat{n}}}{a_{1/2}}, \quad  
  a_{1/2} = \dfrac{a_R + a_L}{2},
  \label{eq:Mach_number}
\end{equation}
where $a_L$ and $a_R$ are defined in \eqref{eq:vel_sound}.
Then, $\widetilde{p}$ is given by
\begin{equation}
\widetilde{p} = \mathscr{P}^{+}_{1}(M_L)p_L + \mathscr{P}^{-}_{1}(M_R)p_R + p_{u}, \label{eq:pIntAusm}
\end{equation}
where $\mathscr{P}^{\pm}_{1}(\cdot)$ is a polynomial defined as \cite{LIOU1996364}:
\begin{equation}
\mathscr{P}^{\pm}_{1}(M) =  \begin{cases}
\frac{1}{M}\mathcal{M}^{\pm}_{1} &\quad\text{if}\,\,\,|M|\geq 1, \\
\\ \mathcal{M}^{\pm}_{2}\left[ \left(\pm 2-M \right)\mp 3 M \mathcal{M}^{\mp}_{2}\right] &\quad\text{otherwise},
\end{cases}
\label{eq:PressurePolynomialAUSM}
\end{equation}
with 
\begin{align} \label{eq:MachPolynomialAUSM}
    \mathcal{M}^{\pm}_{1} = \frac{1}{2}(M \pm |M|), \quad \mathcal{M}^{\pm}_{2}= \pm \frac{1}{4}(M \pm 1)^2.
\end{align}
In eq.~\eqref{eq:pIntAusm}, $p_u$ is a diffusion term introduced to damp the pressure oscillations generated in the limit $M \rightarrow 0$. See \cite{LIOU2006137} for details.






To  compute the mass flow rate $\dot{m}$, we need to define
\begin{equation}
    \widetilde{M} = \mathcal{\mathscr{P}^{+}}_{2}(M_L) +
    \mathcal{\mathscr{P}^{-}}_{2}(M_R) + M_p,\label{eq:MachInterfaceAUSM}
\end{equation}
where $\mathscr{P}^{\pm}_{2}(\cdot)$ is a polynomial defined as: 
\begin{equation}
\mathscr{P}^{\pm}_{2}(M) = \begin{cases}
\mathcal{M}^{\pm}_{1} &\quad\text{if}\,\,\,|M|\geq 1, \\
\mathcal{M}^{\pm}_{2}(1\mp\mathcal{M}^{\mp}_
{2}) &\quad\text{otherwise}.
\end{cases}
\end{equation}
Just like $p_u$ in \eqref{eq:pIntAusm}, $M_p$ in \eqref{eq:MachInterfaceAUSM} is a diffusion term. See \cite{LIOU2006137} for details. Finally, $\dot{m}$ is given by 
\begin{equation}
   \dot{m} = 
   \begin{cases}
    a_{1/2}\widetilde{M} \rho_L &\quad\text{if}\,\,\widetilde{M}>0, \\
    a_{1/2}\widetilde{M} \rho_R &\quad\text{otherwise}.
   \end{cases}
    \label{eq:mDotAUSM}
\end{equation}



\subsubsection{The HLLC-AUSM method}
The HLLC-AUSM method \cite{LIOU2006137,KitKeiiShima} combines the solution of the HLLC and AUSM approximate Riemann solvers to create an extension of the HLLC method suited 
for low Mach number flows. 

The numerical flux $\mathbf{f}_j$ is computed as follows: 

\begin{equation}
    \mathbf{f}_{j} = \frac{\dot{m} + |\dot{m}|}{2} \mathbf{\Phi'_L} + \frac{\dot{m} - |\dot{m}|}{2} \mathbf{\Phi'_R} + \widetilde{p} \mathbf{N}, \label{eq:fluxAUSM+-up2}
\end{equation}
where
\begin{equation*}
 \mathbf{\Phi'_L} =  \mathbf{\Phi_L} + \begin{pmatrix} 
 0 \\ \mathbf{0} \\ \dfrac{{S_L}\left(p^* - p_L\right)}{\rho_L \left(S_L - \mathbf{u}_L \cdot \mathbf{\widehat{n}}\right)} \\  \end{pmatrix}, \quad
 \mathbf{\Phi'_R} =  \mathbf{\Phi_R} + \begin{pmatrix} 0 \\ \mathbf{0} \\ \dfrac{{S_R}\left(p^* - p_R\right)}{\rho_R \left(S_R - \mathbf{u}_R \cdot \mathbf{\widehat{n}}\right)} \\  \end{pmatrix}. 
\end{equation*}

Then the mass flow rate $\dot{m}$ is given by:
\begin{equation}
   \dot{m} = 
   \begin{cases}
    \rho_L \mathbf{u}_L \cdot \mathbf{\widehat{n}} + S_L \left( \rho_L^* - \rho_L\right) \quad \text{if }S^* > 0,\\
    \rho_R \mathbf{u}_R \cdot \mathbf{\widehat{n}} + S_R \left( \rho_R^* - \rho_R\right) \quad \text{otherwise},
   \end{cases}
    \label{eq:mDotHLLCAUSM}
\end{equation}
where 
\begin{equation}
\rho_\alpha^* = \dfrac{S_\alpha - \mathbf{u}_\alpha \cdot \mathbf{\widehat{n}}}{S_\alpha - S^*}\rho_\alpha, \quad \alpha = R, L.
    \label{eq:rhostar}
\end{equation}

\section{Time discretization} \label{sec:time_disc}

While space discretization of problem \eqref{eq:discreteSystem_visc} is rather involved, time 
discretization is rather simple. 
This is the reason why we chose to present space discretization first.

Let us start by noting that we can conveniently rewrite the space-discrete version of system \eqref{eq:discreteSystemII} as:
\begin{equation}\label{eq:space_disc_pb}
   \frac{\partial{\mathbf{q}_i}}{\partial t} = \cL_i(\mathbf{q}),
\end{equation}
with 
\begin{equation*}
\cL_i(\mathbf{q}) = - \sum_j \widetilde{\bS}_{0,j} + \mathbf{\Gamma} \sum_j \nabla \mathbf{D}^{}_j -  \sum_j \mathbf{f}_{j} {A}_j.
\end{equation*}
To discretize \eqref{eq:space_disc_pb} in time, 
we consider a time step $\Delta t \in R^+$.  Let $t^n = n \Delta t$ with $n = 0, 1, \dots, N_T$ and
$t_f = N_T \Delta t$.  We denote with $f^n$ be the approximation of generic variable $f$ at time $t^n$. 

We adopt a fully explicit fourth-order Runge-Kutta scheme \cite{quarteroni2009Book}
that applied to eq.~\eqref{eq:space_disc_pb} reads: find $\mathbf{q}^{n+1}_i$ such that
\begin{equation}
  \mathbf{q}^{n+1}_i = \mathbf{q}^{n}_i + \frac{\Delta t}{6} \left( \bK_1 + 2 \bK_2 + 2 \bK_3 + \bK_4 \right), 
\end{equation}
where
\begin{align}
    &\bK_1 = \cL_i(t^n,\mathbf{q}^n), \quad \bK_2 = \cL_i(t^n + \frac{\Delta t}{2},\mathbf{q}^n + \Delta t \frac{\bK_1}{2}), \cl
    &\bK_3 = \cL_i(t^n + \frac{\Delta t}{2},\mathbf{q}^n + \Delta t \frac{\bK_2}{2}), \quad
    \bK_4 = \cL_i(t^n + \Delta t,\mathbf{q}^n + \Delta t \bK_3). \el
\end{align}
We chose this scheme because it introduces little numerical dissipation. Lower order methods, like, e.g., BDF2, would lead to over-diffusion. We note that this is not the case for a pressure-based solver. See, e.g., 
\cite{GQR_OF_clima} for numerical results
obtained with a pressure-based solver and BDF1 that do not display over-diffusion. 


\begin{rem}
In  Sec.~\ref{sec:space_disc}, we have described how one could build a well-balanced scheme in theory. However, in order to mitigate the numerical error \cite{KAPPELI2014199, bottaKlein2004}, the local hydrostatic profiles, $\rho_{0,i}(z)$ and $p_{0,i}(z)$, defined in \eqref{eq:rhoRec}, \eqref{eq:pGamma} and used in \eqref{eq:source}, \eqref{eq:pCompI}-\eqref{eq:rhoCompI}, \eqref{eq:pert1}-\eqref{eq:pert2} are updated at each time step: they become $p_{0,i}(z, t^n)$ and $\rho_{0,i}(z, t^n)$ computed from \eqref{eq:rhoRec}, \eqref{eq:pGamma} where $\rho_0(z_i)$ is replaced by the cell average density at the current time instant $\rho_0(z_i, t^n)$. 
    Notice that ${p}_{i}'({x}_{i}, t^n) = {\rho}_{i}'({x}_{i}, t^n) = 0 $, i.e., the equilibrium reconstructions 
equal the total cell averages 
\cite{KAPPELI2014199}. 
\end{rem}

\section{Numerical Results}\label{sec:num_res}

The goal of this section is to 
compare the accuracy of the different methods for the computation of the
numerical flux presented
in Sec.~\ref{sec:flux}. For this, 
we consider two classical  benchmarks: the smooth rising thermal bubble \cite{restelliPHD2007,ResGir2009}
and the density current \cite{ahmadLindeman2007,strakaWilhelmson1993}.
Both tests involve a perturbation of a neutrally
stratified atmosphere with uniform background potential temperature over a flat terrain. Therefore, before reporting the results for
the two benchmarks, in Sec.~\ref{sec:hydro}  we show that an unperturbed
stratified atmosphere with uniform background potential temperature over a flat terrain remains unchanged up to a certain tolerance. This is important for illustrating the well-balanced property of our solver. Then, our results for the smooth rising thermal bubble benchmark are presented in Sec.~\ref{sec:bubble}, while Sec.~\ref{sec:DC} reports the results for the density current test. 

We would like to point out that neither the rising bubble nor
the density current benchmark has an exact solution. Hence, one can
only have a relative comparison with other numerical data available
in the literature. 

All the simulations in this section have been carried out with GEA \cite{GEA}.

\subsection{Hydrostatic atmosphere}\label{sec:hydro}

We consider an initial resting atmosphere over a flat terrain. A well-balanced scheme is expected to maintain the atmosphere still for a long time interval with a reasonable accuracy. 
The computational domain in the $xz$-plane is $\Omega$ = [0, 16000] m$\times$ [0, 800] m$^2$. In this domain, the hydrostatic atmosphere, initially
at rest, is free to evolve until $t$ = 25 days \cite{bottaKlein2004, GQR_OF_clima, marrasNazarovGiraldo2015}. 
We impose a no-flux boundary condition at  all the boundaries.
The initial potential temperature is $\theta^0 = 300$ K.

\begin{figure}[htb!]
\centering
\begin{overpic}[abs,unit=1mm,scale=.2]{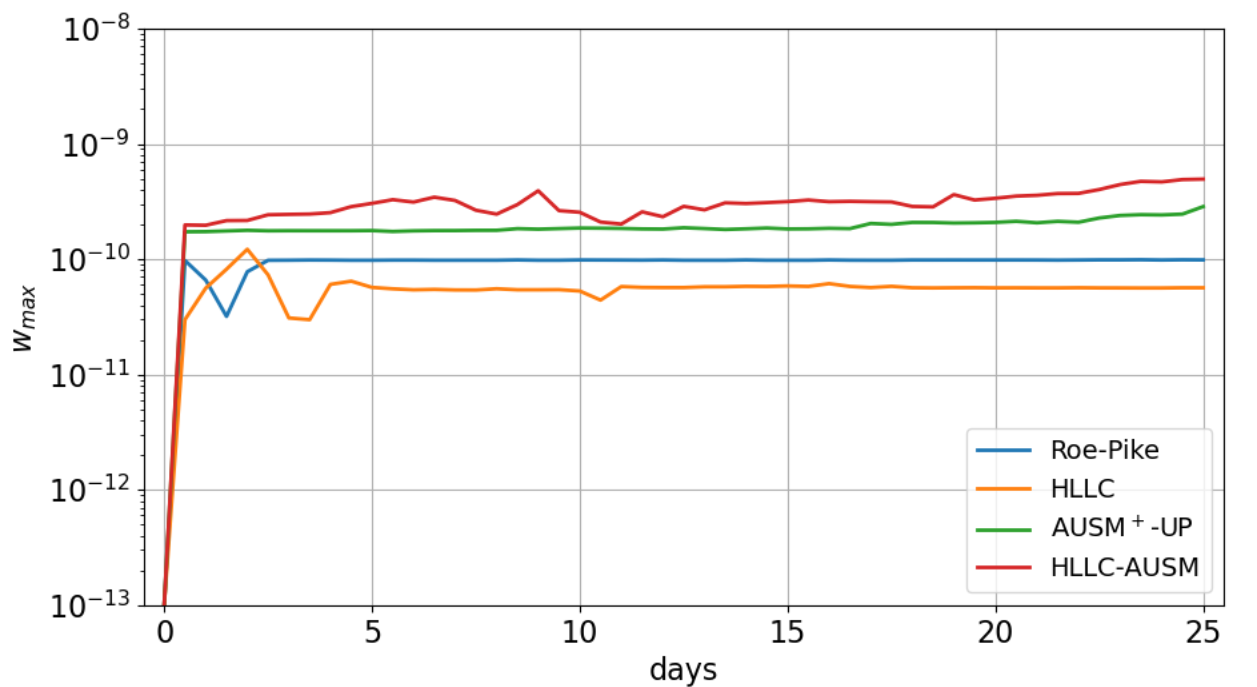}
\end{overpic}
\caption{Hydrostatic atmosphere: time evolution of the maximal vertical velocity
$w_{max}$ for all the methods for the computation of the numerical flux under consideration.
}
\label{fig:FLAT_TERRAIN_uyMax}
\end{figure}

We consider
a uniform mesh with mesh size $h = \Delta x = \Delta z = 250$ m \cite{bottaKlein2004, GQR_OF_clima}
and we set the time step to $\Delta t$ = 0.1 s. Fig.~\ref{fig:FLAT_TERRAIN_uyMax} shows the time evolution of the maximal vertical velocity $w_{max}$. We see 
that in the ``worst'' case (HLLC-AUSM method) $w_{max}$ does not exceed $1e-9$ m/s over the
25 day period. All the other methods keep the value of 
$w_{max}$ even lower. 

We conclude that all methods to compute the numerical flux 
considered in this paper preserve
the hydrostatic equilibrium with reasonably good accuracy.

\subsection{Smooth rising thermal bubble}\label{sec:bubble}

For this benchmark, the computational domain in the $xz$-plane is $\Omega= [0,1000]\times[0,1000]$ m$^2$ and the time interval of interest is (0, 600] s. The initial potential temperature profile is 
\begin{equation}
    \theta^0 = 300 +\frac{0.5}{2}\left[ 1 +\cos\left(\frac{\pi r}{r_c} \right) \right] ~ \textrm{if $r\leq r_c=250~\mathrm{m}$}, \quad\theta^0 = 300
~ \textrm{otherwise},
\label{warmBubble}
\end{equation}
where $r = \sqrt[]{(x-x_{c})^{2} + (z-z_{c})^{2}}$, $(x_c,z_c) = (5000,2000)~\mathrm{m}$ is the radius of the circular perturbation. The local (i.e., in each cell) initial density is given by eq.~\eqref{eq:rhoRec}. 
The initial velocity field is zero everywhere. 
The initial total energy is given by: 

\begin{equation}
    e^0 = U^0 + \Phi 
    \label{eq:InitRTB_Energy}
\end{equation}
where $\Phi = gz$ and $U^0 = c_{v} T^0$, with $T^0$ that can be computed from $\theta^0$ in \eqref{warmBubble}.
No-flux boundary conditions are
imposed on all walls. 

We consider two meshes with uniform resolution: $h = \Delta x = \Delta z=2.5, 5$ m. We set $\Delta t = 0.05$ s. 
 Furthermore, following \cite{GEA_GIR_QUA} in \eqref{eq:momentum_stab}-\eqref{eq:energy_stab} we set $\mu_a=0.15$ m$^2$/s and $Pr=1$.
We compare our numerical results with the results reported in \cite{ResGir2009,GEA_GIR_QUA}. 
The results in \cite{ResGir2009} are obtained with a density-based approach developed from a Godunov-type scheme, similar to the approach used in this paper, but 
Discontinuous Galerkin and Spectral
Elements methods are used for space discretization. For this reason, we considered also the results from \cite{GEA_GIR_QUA}, which were obtained 
with a different approach (i.e., pressure-based)
but with the same space discretization
method (i.e., a Finite Volume method in GEA \cite{GEA}).

\begin{figure}[htb!]
\centering
           \begin{overpic}[width=0.43\textwidth]{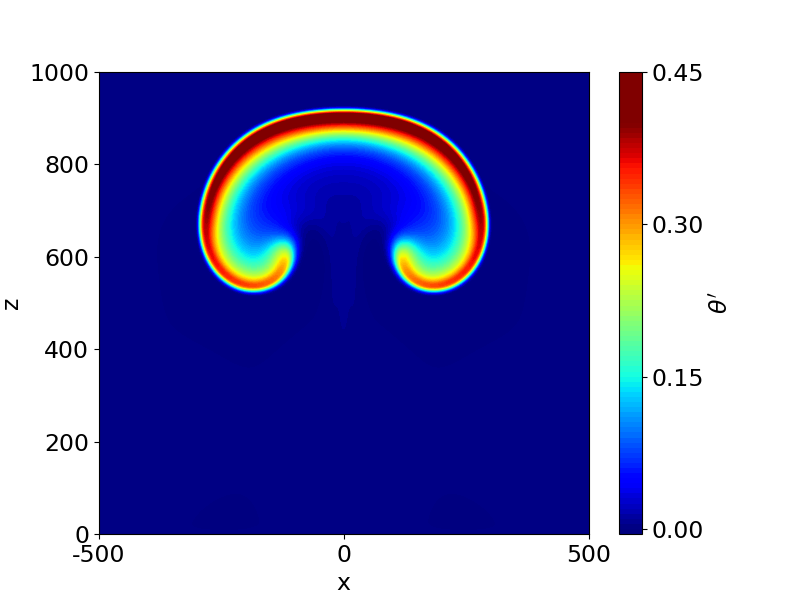}
    \put(63,28){\textcolor{white}{Roe-Pike}}
    \end{overpic} 
    \begin{overpic}[width=0.425\textwidth]{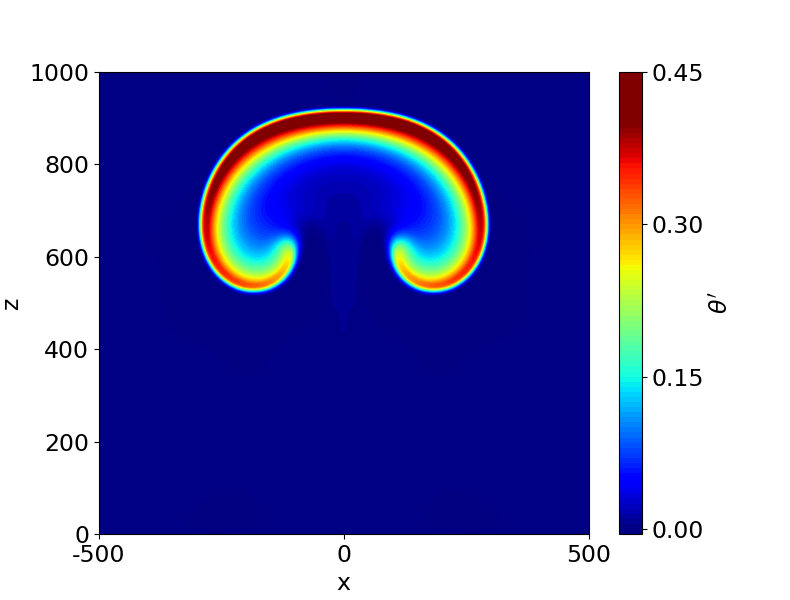}
     \put(71,28){\textcolor{white}{HLLC}}
    \end{overpic} \\ \vspace{0.5cm}
    \begin{overpic}[width=0.43\textwidth]{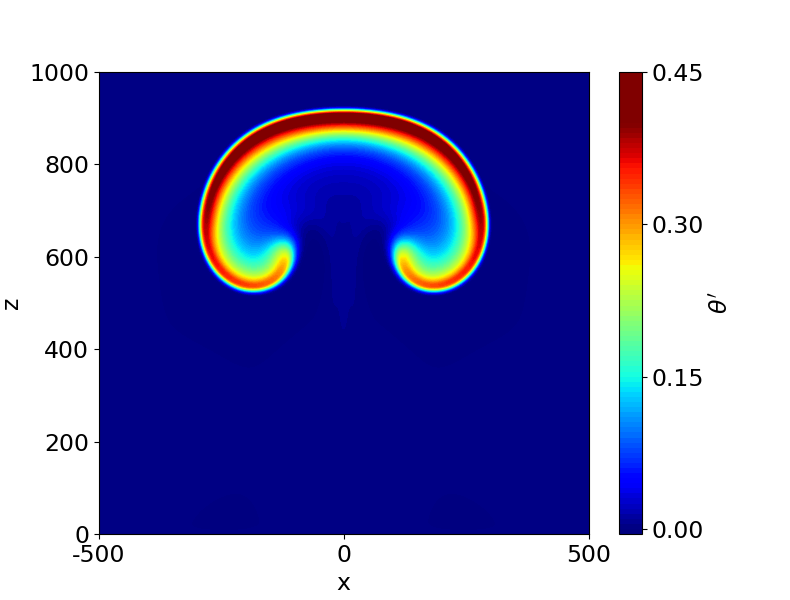}
    \put(60,28){\textcolor{white}{AUSM$^+$-up}}
    \end{overpic} 
    \begin{overpic}[width=0.43\textwidth]{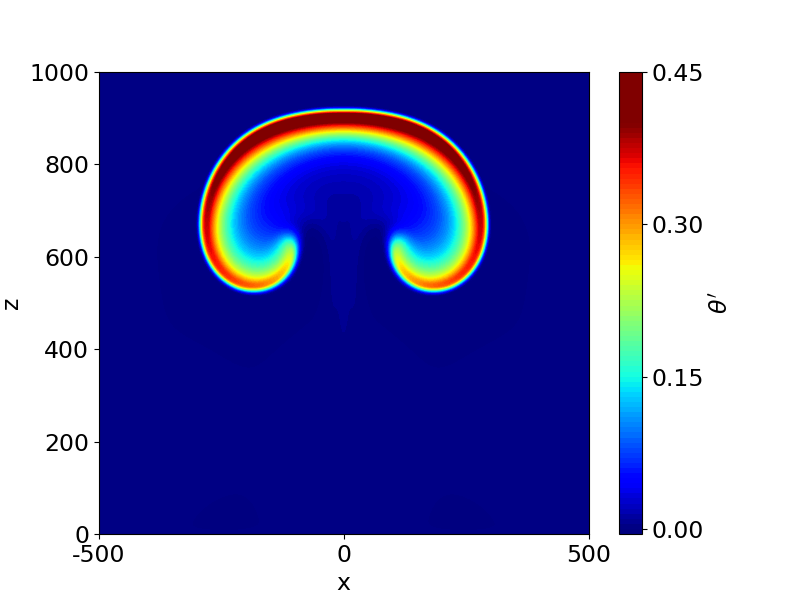}
       \put(55,28){\textcolor{white}{HLLC-AUSM}}
    \end{overpic}
    \caption{Rising thermal bubble: perturbation of the potential temperature computed
at $t = 600$ s with mesh $h= 2.5$ m and 
the Roe-Pike (top left), HLLC (top right), AUSM$^{+}$-up (bottom letft), and HLLC-AUSM (bottom right) methods.
}
    \label{fig:snapRTB}
\end{figure}
 
\begin{figure}[htb!]
    \centering
    \begin{overpic}[width=0.9\textwidth]{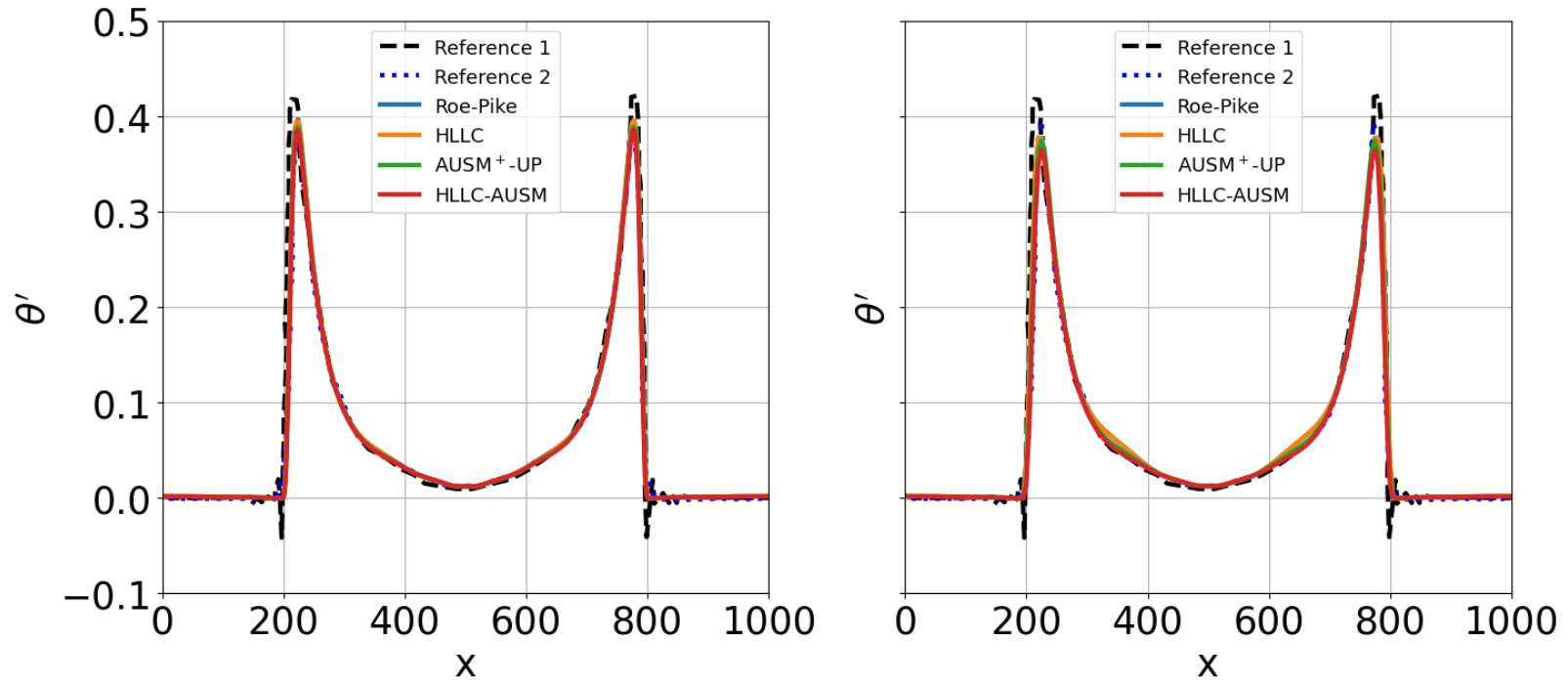}
    \put(98,182){\textcolor{black}{\large{$h=2.5$ m}}}
    \put(300,182){\textcolor{black}{\large{$h=5$ m}}}
    \end{overpic} 
    \caption{Rising thermal bubble: profile of the potential temperature perturbation along $z=700$ m at $t=600$ s given by the different 
    methods to compute the numerical flux with 
    mesh $h=2.5$ m (left) and mesh $h=5$ m (right). Reference 1 is taken from 
    \cite{ResGir2009}, while Reference 2 is from \cite{GEA_GIR_QUA}.
}
    \label{fig:thetaRTB2}
\end{figure}

Fig.~\ref{fig:snapRTB} shows the potential temperature perturbation computed at $t=600$ s with mesh $h=2.5$ m and all the methods for the computation of
numerical flux in Sec.~\ref{sec:flux}. 
We notice no visible difference in the solutions given by the different methods. Moreover, we have
very good qualitative agreement with the results in \cite{ResGir2009,GEA_GIR_QUA}. 
For an easier comparison with the results in \cite{ResGir2009,GEA_GIR_QUA},
Fig. \ref{fig:thetaRTB2} depicts the profile of the potential temperature perturbation
along $z = 700$ m at $t=600$ s. The reference values from \cite{ResGir2009,GEA_GIR_QUA} correspond to mesh size $h = 5$ m, while  we show the results for $h = 2.5, 5$ m.
We see that all the curves are rather close to each other, with the following minor exceptions: around the maxima for
both $h = 2.5, 5$ m and 
around $x = 300$
and $x = 700$ for $h = 5$ m. In addition, 
our approach does not exhibit
the same oscillations, supposedly numerical, around $x = 200$
and $x = 800$ as in \cite{ResGir2009}.

For a more quantitative comparison, Table \ref{tab:RTB2_uv} reports the maximum and the minimum value of the velocity components. 
Overall, 
the best agreement with the data from \cite{ResGir2009,GEA_GIR_QUA} is for the values
given by the HLLC-AUSM method. 
The AUSM$^{+}$-up scheme
is slightly more dissipative, while the ROE-Pike and
HLLC methods dissipate even more. 


\begin{table}[htb!]
\begin{center}
\begin{tabular}{ | c | c |  c | c | c | c | }
\hline
Method & $h$ (m) & $u_{min}$ (m/s) & $u_{max}$ (m/s) & $w_{min}$ (m/s) & $w_{max}$ (m/s)  \\
 \hline
Ref.~\cite{ResGir2009}  & 5 & -2.16 & 2.16 & -1.97 & 2.75 \\
 \hline
 Roe-Pike & 5 & -1.65 & 1.65 & -1.60 & 2.47 \\
Roe-Pike &  2.5 & -1.80 & 1.80 & -1.68 & 2.50 \\
 HLLC & 5  & -1.62 & 1.62 & -1.60 & 2.46 \\

 HLLC &  2.5 & -1.80 & 1.80 & -1.68 & 2.50 \\
 
 AUSM$^+$-up & 5  & -1.75 & 1.75 & -1.65 & 2.50 \\

 AUSM$^+$-up & 2.5 & -1.87 & 1.87 & -1.70 & 2.50 \\
 
 HLLC-AUSM  & 5  & -1.85 & 1.85 & -1.69 & 2.48 \\

 HLLC-AUSM & 2.5 & -1.92 & 1.92 & -1.71 & 2.51 \\ 
 \hline
\end{tabular}
\caption{Rising thermal bubble: maximum and minimum values of the horizontal component $u$ and vertical component $w$ of the velocity at $t=600$ s. 
}\label{tab:RTB2_uv}
\end{center}
\end{table}

\subsection{Density current}\label{sec:DC}

The computational domain of this benchmark is $\Omega=[0,25600]\times[0,6400]~\mathrm{m}^2$ in the $xz$-plane and 
the time interval of interest is $(0,900]$ s. 
The initial potential temperature profile is
\begin{equation}
\theta^0 = 300 - \frac{15}{2}\left[  1 + \cos(\pi r)\right] ~ \textrm{if $r\leq 1$},\quad\theta^0 = 300
~ \textrm{otherwise},
\label{dcEqn1}
\end{equation}
where $r = \sqrt[]{\left(\frac{x-x_{c}}{x_r}\right)^{2} + \left(\frac{z-z_{c}}{z_r}\right)^{2}}$, with $(x_r,z_r)=(4000, 2000)~{\rm m}$ and $(x_c,z_c) = (0,3000)~\mathrm{m}$. 
The local initial density field is given by eq.~\eqref{eq:rhoRec}. The initial velocity field is zero everywhere. The initial total energy is given by eq. \eqref{eq:InitRTB_Energy}. No-flux boundary conditions are imposed on all walls.

We consider four meshes with uniform resolution: $h =\Delta x$ = $\Delta z  = 200, 100, 50, 25$ m. We set $\Delta t=0.05$ s. Furthermore, in \eqref{eq:momentum_stab}-\eqref{eq:energy_stab} we set $\mu_a=75$ m$^2$/s and $Pr=1$ as done in \cite{GQR_OF_clima,strakaWilhelmson1993,ahmadLindeman2007}. 





\begin{figure}[htb]
\centering
    \begin{overpic}[width=0.99\textwidth,grid = false]{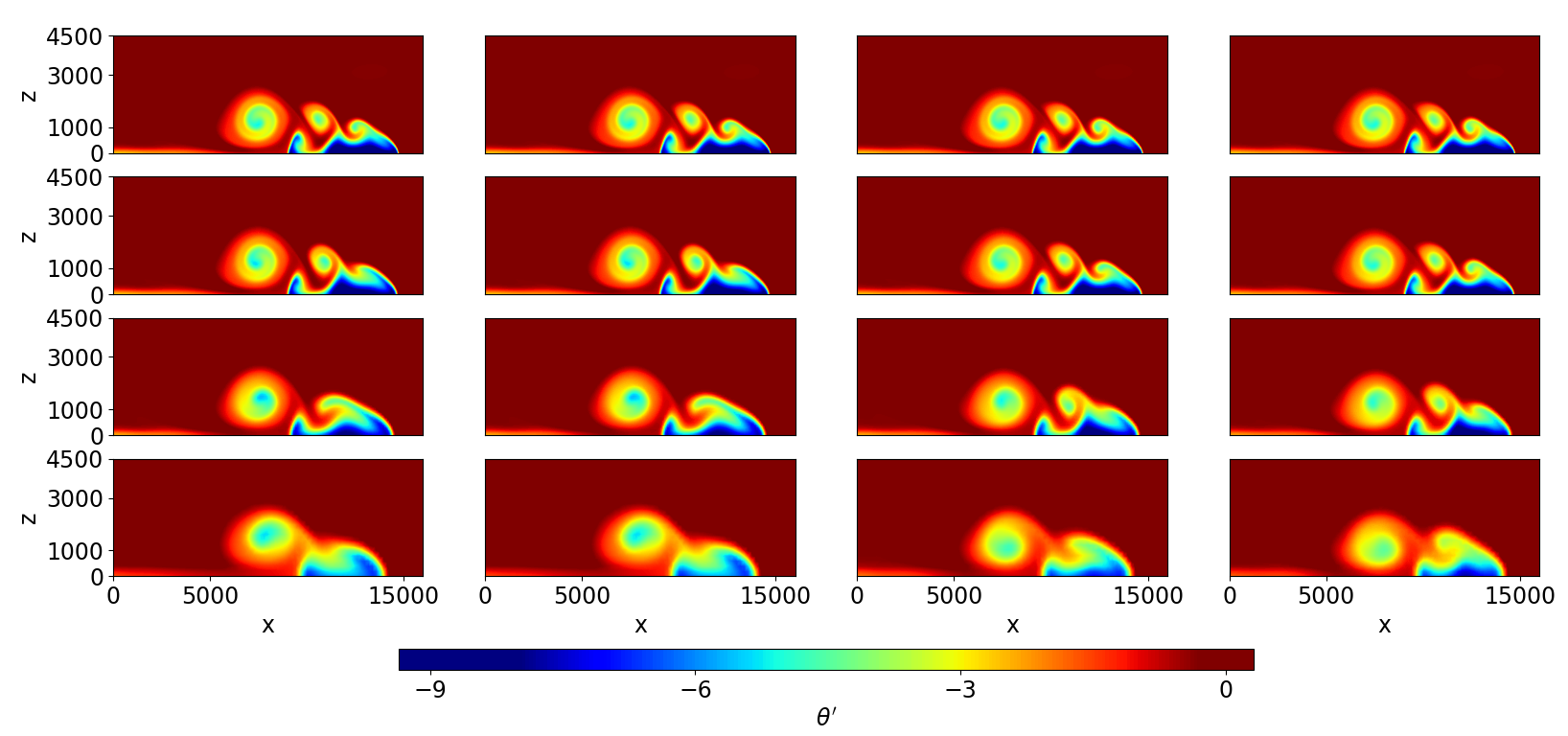}
    \put(55,215){\textcolor{black}{Roe-Pike}}
    \put(175,215){\textcolor{black}{HLLC}}
    \put(270,215){\textcolor{black}{AUSM$^+$-up}}
    \put(370,215){\textcolor{black}{HLLC-AUSM}}
    \put(40,73){\textcolor{white}{\tiny{$h=200$ m}}}
    \put(40,115){\textcolor{white}{\tiny{$h=100$ m}}}
    \put(40,155){\textcolor{white}{\tiny{$h=50$ m}}}
    \put(40,199){\textcolor{white}{\tiny{$h=25$ m}}}%
    \put(150,73){\textcolor{white}{\tiny{$h=200$ m}}}
    \put(150,115){\textcolor{white}{\tiny{$h=100$ m}}}
    \put(150,155){\textcolor{white}{\tiny{$h=50$ m}}}
    \put(150,199){\textcolor{white}{\tiny{$h=25$ m}}} %
    \put(370,73){\textcolor{white}{\tiny{$h=200$ m}}}
    \put(370,115){\textcolor{white}{\tiny{$h=100$ m}}}
    \put(370,155){\textcolor{white}{\tiny{$h=50$ m}}}
    \put(370,199){\textcolor{white}{\tiny{$h=25$ m}}}
    \put(260,73){\textcolor{white}{\tiny{$h=200$ m}}}
    \put(260,115){\textcolor{white}{\tiny{$h=100$ m}}}
    \put(260,155){\textcolor{white}{\tiny{$h=50$ m}}}
    \put(260,199){\textcolor{white}{\tiny{$h=25$ m}}}
    \end{overpic}
    \caption{Density current: potential temperature fluctuation $\theta'$  given by all the numerical methods for the flux computation at $t = 900$ s with meshes $h = 25$ m (first row), $h = 50$ m (second row), $h = 100$ m (third row), $h = 200$ m (fourth row). 
    }
    \label{fig:DC_flux_comparison}
\end{figure}
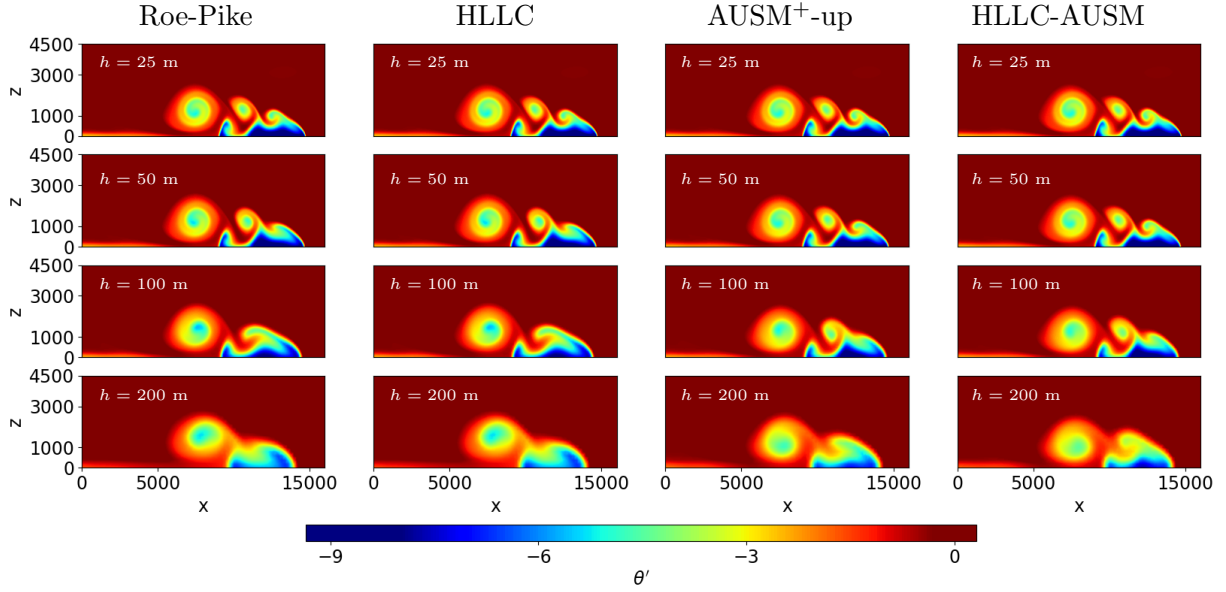

Fig.~\ref{fig:DC_flux_comparison} illustrates the perturbation of potential temperature $\theta'$ at $t = 900$ s computed with all the meshes and all the methods for the computation of the numerical flux in Sec.~\ref{sec:flux}.
With the finest mesh ($h = 25$ m), all the methods are able to 
capture a clear three-rotor structure, which is in good agreement with the 
results reported in the literature for the same resolution. 
See, e.g., \cite{ahmadLindeman2007,giraldo_2008,GQR_OF_clima,marrasEtAl2013a,marrasNazarovGiraldo2015,strakaWilhelmson1993}. 
With mesh $h = 50$ m, the three-rotor structure is still well captured
by the AUSM$^+$-up and HLLC-AUSM  methods, while the Roe-Pike and HLLC
methods dampen the smallest recirculation. %
With mesh $h=100$ m, the smallest recirculation is significantly damped 
by the AUSM$^+$-up and HLLC-AUSM methods too. With the same mesh, 
the Roe-Pike and HLLC solutions have the main rotor and a prolonged
recirculation resulting from the merging of the smaller rotors. 
Finally, we see that mesh $h=200$ m is too coarse and no method
is able to provide an accurate solution. 
From now on, we will consider only meshes $h = 25, 50$ m.

\begin{figure}[htb!]
    \centering
    \begin{overpic}[width=0.45\textwidth]{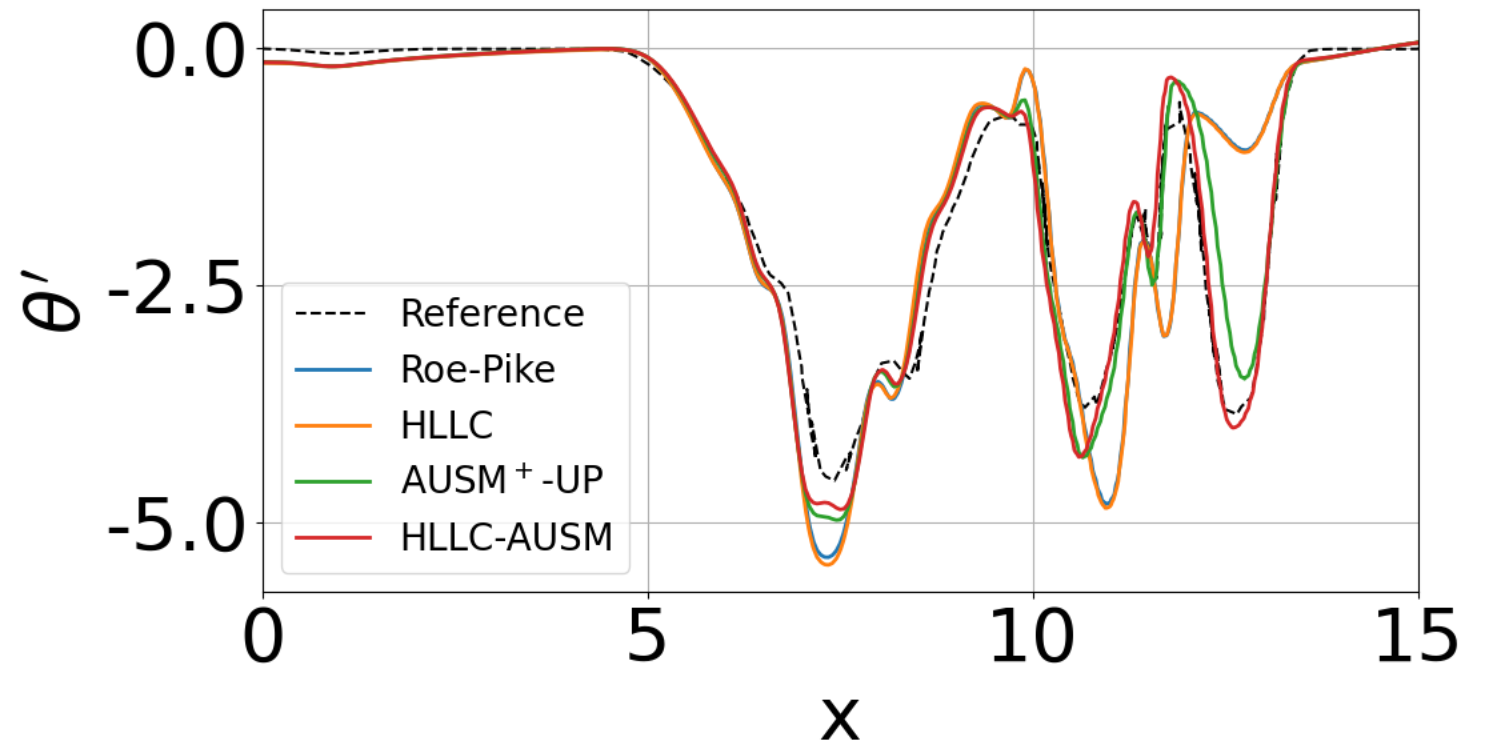}
    \put(90,109){\textcolor{black}{$h=50$ m}}
    \end{overpic}  
     \begin{overpic}[width=0.45\textwidth]{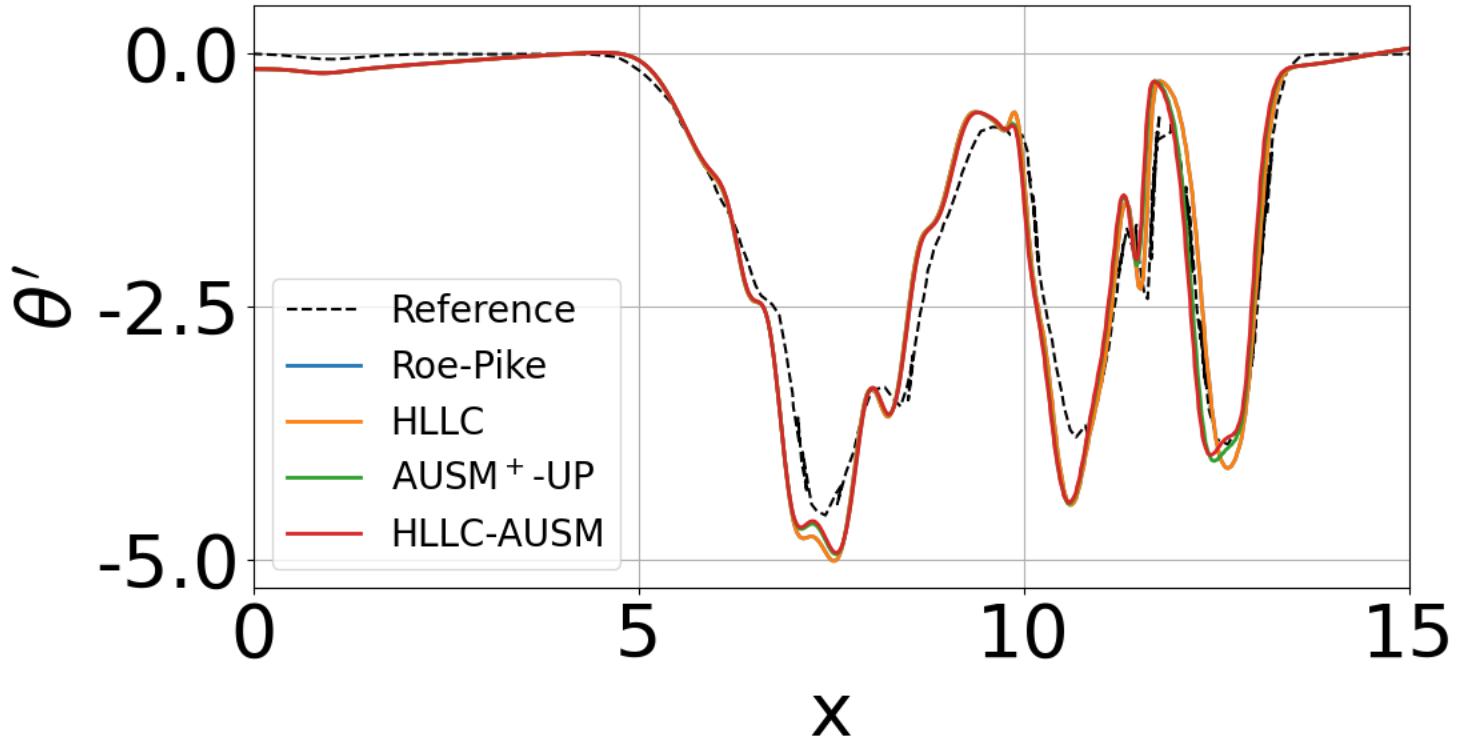}
     \put(90,109){\textcolor{black}{$h=25$ m}}
    \end{overpic}
    
    \caption{Density current: potential temperature perturbation $\theta'$ along $z = 1200$ m  at $t = 900$ s given by all the methods for the computation of the numerical flux
    with mesh $h = 50$ m (left) and $h = 25$ m (right).
    The Reference data are taken from \cite{giraldo_2008} and refer to resolution 25 m. 
}
\label{fig:DC_FLUXES_COMP_1200}
\end{figure}

Fig.~\ref{fig:DC_FLUXES_COMP_1200} compares the potential temperature perturbation along $z = 1200$ m at $t = 900$ s 
computed with meshes $h = 50, 25$ m and the data in \cite{giraldo_2008}, 
which refer to resolution 25 m. We recall that the data in \cite{giraldo_2008} are obtained
with a spectral element method and a discontinuous Galerkin method. For this benchmark, these two methods give results 
so close that the corresponding curves appear
superimposed.
For this reason, we label them just as Reference
in Fig.~\ref{fig:DC_FLUXES_COMP_1200}. 
Let us first comment on the results in Fig. \ref{fig:DC_FLUXES_COMP_1200} obtained for 
mesh $h = 50$ m (left panel). With the HLLC-AUSM method, we get results that are in good 
agreement with the reference values, with the exception of the negative peaks associated with the two larger rotors. The same mismatch at these negative 
peaks is observed also for the AUSM$^+$-up method, which additionally
gives a smaller (in absolute value) negative peak for
the smallest recirculation.
The curves related to the HLLC and Roe-Pike methods, which are 
almost superimposed, show larger differences at the negative peaks and 
are off phase with respect to the Reference.
Looking at the results for mesh $h = 25$ m (right panel in Fig. \ref{fig:DC_FLUXES_COMP_1200}), we see that the curves obtained
with the different methods are almost superimposed. 
With the exception of the amplitude of the negative peaks,
they are in very good agreement with the reference values. 
This is remarkable if we consider that the results in \cite{giraldo_2008}
are obtained with high-order methods, while we use a
second order accurate finite volume method.  

\vskip .5 cm

The front location for this benchmark is defined as the location on the ground where $\theta'=-1$ K. Table \ref{tab:FrontLoc} reports the front location  
at $t = 900$ s computed with meshes $h = 25, 50$ m. In the table, 
we report also the range of
mesh sizes (from 25 m
to 200 m) and front locations obtained with 14 different methodologies
in \cite{strakaWilhelmson1993}. We observe that
in all the cases our results fall well within the values
reported in~\cite{strakaWilhelmson1993}. 
Additionally, we note that the front locations
obtained with mesh $h = 25$ m are all within 20 m of each other, in a domain 
that is 25.6 Km long. This increases to 65 m with mesh $h = 50$ m, 
which is still pretty good given the size of the domain.

\begin{table}[htb!]
\begin{center}
\begin{tabular}{ | c | c |  c | }
\hline
Method & Resolution [m] & Front Location [m]  \\

 \hline
Roe-Pike & 50 & 14724 \\
\hline
Roe-Pike & 25 & 14780 \\
 \hline
 HLLC & 50  & 14720 \\
 \hline
 HLLC & 25  & 14780  \\
 \hline
 AUSM$^+$-up  & 50  & 14885  \\
 \hline
 AUSM$^+$-up  & 25  & 14790  \\
\hline
 HLLC-AUSM  & 50  & 14765  \\
 \hline
 HLLC-AUSM  & 25  & 14800  \\
 \hline
    Ref.~\cite{strakaWilhelmson1993} & (25, 200) & (14533,17070)
 \\
  \hline
\end{tabular}
\caption{Density current: front location at $t = 900$ s obtained with the different methods and meshes $h = 50, 25$ m.
For reference \cite{strakaWilhelmson1993}, 
we reported the range of
mesh sizes and front location values obtained with 
different methods.}
\label{tab:FrontLoc}
\end{center}
\end{table}

\section{Concluding remarks}\label{sec:conclusion}

We developed, and implemented in a Finite Volume environment, 
a well-balanced density-based solver for the numerical 
simulation of non-hydrostatic atmospheric flows. To approximate the solution
of the Riemann problem, we considered four methods: 
Roe-Pike, HLLC, AUSM$^+$-up, and HLLC-AUSM. 
We assessed our density-based approach and compared 
the accuracy of these four methods through two well-known
benchmarks: the smooth rising thermal bubble and the density current.

We found that the solutions given by the different approximated Riemann solvers
differ noticeably when using coarser meshes. Specifically, 
unless the mesh is very fine, the Roe-Pike and HLLC methods give over-diffusive solutions,
while both the AUSM$^+$-up and the HLLC-AUSM methods are less dissipative
and thus allow for the use of coarser meshes. In particular, the 
HLLC-AUSM method is the one that gives the best comparison with the data
available in the literature, even with coarser meshes. 
The differences in the solutions given by the approximated Riemann solvers
become less evident as the mesh gets refined.

\section*{Aknowledgements}
We acknowledge the support provided by PRIN “FaReX - Full and Reduced order modelling of
coupled systems: focus on non-matching methods and automatic learning” project, PNRR NGE iNEST “Interconnected Nord-Est Innovation Ecosystem” project, INdAM-GNCS 2019–2020 projects
and PON “Research and Innovation on Green related issues” FSE REACT-EU 2021 project. This
work was also partially supported by the U.S. National Science Foundation through Grant No.
DMS-1953535 (PI A. Quaini).



\bibliographystyle{plain}
\bibliography{bibliography.bib}

\end{document}